\documentclass[11pt]{article}
\usepackage{amssymb,amsthm,amsmath}
\usepackage{cite}

\def\qed{\hfill $\Box$}

\newcommand\pf{\smallbreak\noindent \texttt{Proof}. }

\begin{document}

\newtheorem{thm}{Theorem}[section]
\newtheorem{prop}[thm]{Proposition}
\newtheorem{lem}[thm]{Lemma}
\newtheorem{cor}[thm]{Corollary}
\newtheorem{ex}[thm]{Example}
\renewcommand{\thefootnote}{*}

\title{\bf Leibniz Rings: some basic and structural results}

\author{\textbf{L.A.~Kurdachenko, O.O.~Pypka}\\
Oles Honchar Dnipro National University, Dnipro, Ukraine\\
{\small e-mail: lkurdachenko@gmail.com, sasha.pypka@gmail.com}\\
\textbf{M.M.~Semko}\\
State Tax University, Irpin, Ukraine\\
{\small e-mail: dr.mykola.semko@gmail.com}}
\date{}

\maketitle

\begin{abstract}
In this paper, we study the fundamental properties of Leibniz rings. Special attention is given to the structure of Leibniz rings whose additive group is ``small''. The results obtained illustrate a significant difference between the classes of Leibniz rings and Lie rings.
\end{abstract}



\thispagestyle{empty}

\section*{Introduction}
Lie rings are among the oldest objects of study among all non-associative rings. Their linear analogue -- Lie algebras -- constitute one of the most thoroughly investigated types of non-associative algebras. The theory of Lie algebras is one of the most developed algebraic theories, rich in important and interesting results, many of which have become classical. In turn, Lie algebras are a special case of Leibniz algebras. More precisely, Lie algebras are exactly the anticommutative Leibniz algebras.

Leibniz algebras first appeared and were studied in the papers~\cite{BA1965,BA1967,BA1971}, in which they were called \textit{D-algebras}. However, at that time these papers did not attract much attention. Only two decades later did real interest in Leibniz algebras arise. This happened thanks to the work of J.-L.~Loday~\cite{LJ1993} (see also \cite[Section~10.6]{LJ1992}), who ``rediscovered'' these algebras and used the term \textit{Leibniz algebras}, since it was Leibniz who discovered and proved the Leibniz rule for the differentiation of functions. The main motivation for introducing Leibniz algebras was the study of periodicity phenomena in algebraic K-theory. Leibniz algebras turned out to be naturally related to several areas, such as differential geometry, homological algebra, classical algebraic topology, algebraic K-theory, loop spaces, noncommutative geometry, and so on. Nowadays, the theory of Leibniz algebras is one of the actively developing areas of modern algebra. It should be noted that in recent years two monographs~\cite{AOR2020,KPS2024} and numerous papers (see, e.g., \cite{CPSeY2019, KiKuPSu2017, KOP2016, KPS2018, KPS2021, KPS2022, KPS2023, KPV2023, KSS2017, KSS2018A, KSeSu2020, KSeSu2023, KSeSu2024, KSeSu2024A, KSeY2021, KSeYa2023, KSeYa2024, KSuSe2018, KSuY2024, P2022, SSY2022}) have been published, presenting various results of this theory.

It should be noted that there is a very significant difference between Lie rings and Lie algebras. The additive groups of Lie algebras, as well as of Leibniz algebras, have a very simple structure: they are either elementary abelian $p$-groups for some prime $p$, or torsion-free divisible abelian groups. At the same time, the structure of the additive groups of arbitrary Lie rings can be considerably more complicated. This explains the fact that the theory of Lie rings is not as well developed as the theory of Lie algebras. As for the theory of Leibniz rings, it is scarcely developed at all -- the number of papers devoted to Leibniz rings can be counted on the fingers. This naturally raises the question of a systematic and consistent development of the theory of Leibniz rings.

A set $L$ with two binary operations $+$ and $[,]$ is called a \textit{Leibniz ring} (more precisely, a \textit{left Leibniz ring}) if it satisfies the following properties:
\begin{enumerate}
\item[\upshape(i)] $L$ is an abelian group under addition;
\item[\upshape(ii)] $[a,b+c]=[a,b]+[a,c]$ and $[a+b,c]=[a,c]+[b,c]$;
\item[\upshape(iii)] $[a,[b,c]]=[[a,b],c]+[b,[a,c]]$
\end{enumerate}
for all $a,b,c\in L$.

Dually, if instead of the last equality $L$ satisfies
$$[a,[b,c]]=[[a,b],c]-[[a,c],b],$$
then $L$ is called a \textit{right Leibniz ring}.

Note that the classes of left Leibniz rings and right Leibniz rings are distinct. The following simple example demonstrates this.

Let $L=\langle a\rangle\oplus\langle b\rangle$ be the elementary abelian group of order $p^{2}$, where $p$ is a prime. Define the operation $[,]$ on $L$ by the following rule:
$$[a,a]=[a,b]=b,\ [b,a]=[b,b]=0.$$
It is straightforward to check that $L$ becomes a left Leibniz ring. However,
$$0=[b,a]=[[a,a],a]\neq[[a,a],a]+[a,[a,a]]=[a,b]=b.$$

Let $R$ be a right Leibniz ring, and define $[\![a,b]\!]=[b,a]$. Then we have:
\begin{align*}
[\![[\![a,b]\!],c]\!]&=[c,[b,a]]=[[c,b],a]-[[c,a],b]\\
&=[\![a,[\![b,c]\!]]\!]-[\![b,[\![a,c]\!]]\!].
\end{align*}
Thus, this substitution transforms a right Leibniz ring into a left Leibniz ring. Similarly, one can construct a transfer from a left Leibniz ring to a right Leibniz ring.

A ring $L$ is called a \textit{symmetric Leibniz ring} if it is both a left and a right Leibniz ring. In what follows, we shall study in detail some fundamental properties of Leibniz rings.

\section{Basic concepts and results on Leibniz rings}
\begin{prop}\label{P1}
Let $L$ be a left Leibniz ring. Then $L$ is a symmetric Leibniz ring if and only if $[b,[a,c]]=-[[a,c],b]$ for all $a,b,c\in L$.
\end{prop}
The assertion is almost obvious, so we omit the proof.

We prefer to work with left Leibniz rings. In this regard, we note the following interesting property of Leibniz rings.

\begin{prop}\label{P2}
Let $L$ be a left Leibniz ring. Then $[[a,b],c]=-[[b,a],c]$ for all $a,b,c\in L$.
\end{prop}
\pf
We have:
$$[a,[b,c]]=[[a,b],c]+[b,[a,c]],$$
and
$$[b,[a,c]]=[[b,a],c]+[a,[b,c]],$$
or equivalently,
$$[a,[b,c]]=[b,[a,c]]-[[b,a],c].$$
It follows that
$$[[a,b],c]+[b,[a,c]]=[b,[a,c]]-[[b,a],c],$$
and hence
$$[[a,b],c]=-[[b,a],c].$$
\qed

The first important examples of Leibniz rings are Lie rings.

\begin{prop}\label{P3}
Let $L$ be a Lie ring. Then $L$ is a Leibniz ring in which $[a,a]=0$ for every $a\in L$. Conversely, if $L$ is a Leibniz ring in which $[a,a]=0$ for every $a\in L$, then $L$ is a Lie ring.
\end{prop}
\pf
Let $L$ be a Lie ring. Then we have
$$[[a,b],c]+[[c,a],b]+[[b,c],a]=0.$$
It follows that
\begin{align*}
[[a,b],c]&=-[[c,a],b]-[[b,c],a]\\
&=[a,[b,c]]-[[c,a],b]\\
&=[a,[b,c]]+[b,[c,a]]\\
&=[a,[b,c]]-[b,[a, c]].
\end{align*}

Conversely, let $L$ be a Leibniz ring in which $[a,a]=0$ for all $a\in L$. For arbitrary elements $a,b\in L$, we have
$$0=[a+b,a+b]=[a,a]+[a,b]+[b,a]+[b,b]=[a,b]+[b,a].$$
It follows that $[a,b]=-[b,a]$. Then
\begin{align*}
0&=[[a,b],c]-[a,[b,c]]+[b,[a,c]]\\
&=[[a,b],c]+[[b,c],a]-[[a,c],b]\\
&=[[a,b],c]+[[b,c],a]+[[c,a],b]
\end{align*}
for all $a,b,c\in L$. Hence $L$ is a Lie ring.
\qed

Of course, Proposition~\ref{P1} shows that every Lie ring is a symmetric Leibniz ring. Note that the class of symmetric Leibniz rings is strictly larger than the class of Lie rings. The following example demonstrates this.

Let $L=\langle a\rangle\oplus\langle b\rangle$ be the elementary abelian group of order $p^{2}$, where $p$ is a prime. Define the operation $[,]$ on $L$ by the following rule:
$$[a,a]=b,\ [b,a]=[b,b]=[a,b]=0.$$
Let
$$x=k_{1}a+t_{1}b,\ y=k_{2}a+t_{2}b,\ z=k_{3}a+t_{3}b,$$
$0\leq k_{1},t_{1},k_{2},t_{2},k_{3},t_{3}<p$. We have
\begin{align*}
[x,y]&=[k_{1}a+t_{1}b,k_{2}a+t_{2}b]\\
&=[k_{1}a,k_{2}a]+[k_{1}a,t_{2}b]+[t_{1}b,k_{2}a]+[t_{1}b,t_{2}b]\\
&=k_{1}k_{2}[a,a]+k_{1}t_{2}[a,b]+t_{1}k_{2}[b,a]+t_{1}t_{2}[b,b]\\
&=k_{1}k_{2}b.
\end{align*}
Similarly,
$$[x,z]=k_{1}k_{3}b,\ [y,z]=k_{2}k_{3}b.$$
It follows that
$$[x,[y,z]]=[[x,y],z]=[y,[x,z]]=[[x,z],y]=0,$$
so $L$ is a symmetric Leibniz ring. However, $[a,a]=b\neq0$, so $L$ is not a Lie ring.

Let $L$ be a Leibniz ring. As usual, a subset $A$ of $L$ is called a \textit{subring} of $L$ if it is closed under both operations $+$ and $[,]$, and is a Leibniz ring with respect to the restrictions of these operations.

If $A,B$ are the subgroups of the additive group of $L$, then $[A,B]$ denotes the subgroup of the additive group of $L$ generated by all elements $[a,b]$, where $a\in A$, $b\in B$. In particular, a subgroup $A$ of the additive group of $L$ is a subring of $L$ if and only if $[A,A]\leq A$.

A subring $A$ of a Leibniz ring $L$ is called a \textit{left} (respectively, \textit{right}) \textit{ideal} of $L$ if $[y,x]\in A$ (respectively, $[x,y]\in A$) for every $x\in A$ and $y\in L$. In other words, $A$ is a left (respectively, right) ideal of $L$ if and only if $[L,A]\leq A$ (respectively, $[A,L]\leq A$).

A subring $A$ of $L$ is called an \textit{ideal} (more precisely, a \textit{two-sided ideal}) of $L$ if it is both a left ideal and a right ideal; that is, $[x,y],[y,x]\in A$ for all $x\in A$, $y\in L$.

If $A$ is an ideal of $L$, we can speak of the \textit{factor-ring} (or \textit{quotient ring}) $L/A$. It is straightforward to verify that this factor-ring is also a Leibniz ring.

A Leibniz ring $L$ is called \textit{abelian} (or \textit{trivial}) if $[a,b]=0$ for all $a,b\in L$. In particular, every abelian Leibniz ring is a Lie ring.

Let $L$ be a Leibniz ring. Define the \textit{lower central series} of $L$
$$L=\gamma_{1}(L)\geq\gamma_{2}(L)\geq\ldots\gamma_{\alpha}(L)\geq\gamma_{\alpha+1}(L)\geq\ldots\gamma_{\delta}(L),$$
by the following rule: $\gamma_{1}(L)=L$, $\gamma_{2}(L)=[L,L]$, and, recursively,
$$\gamma_{\alpha+1}(L)=[L,\gamma_{\alpha}(L)]$$
for all ordinals $\alpha$, and
$$\gamma_{\lambda}(L)=\bigcap_{\mu<\lambda}\gamma_{\mu}(L)$$
for all limit ordinals $\lambda$. The last term $\gamma_{\delta}(L)$ is called the \textit{lower hypocenter} of $L$. We have: $\gamma_{\delta}(L)=[L,\gamma_{\delta}(L)]$.

If $\alpha=k$ is a positive integer, then $\gamma_{k}(L)=[L,[L,[L,[\ldots,L]]]\ldots]$.

We note the following useful properties of subrings and ideals.

\begin{prop}\label{P4}
Let $L$ be a Leibniz ring. Then the following assertions hold:
\begin{enumerate}
\item[\upshape(i)] if $H$ is an ideal of $L$ and $S$ is a subring of $L$, then $H+S$ is a subring of $L$;
\item[\upshape(ii)] if $H$ is an ideal of $L$, then $[H,H]$ is an ideal of $L$;
\item[\upshape(iii)] if $H$ is an ideal of $L$, then $[L,H]$ is a subring of $L$;
\item[\upshape(iv)] if $H$ is an ideal of $L$, then $[H,L]$ is a subring of $L$;
\item[\upshape(v)] if $H$ is an ideal of $L$, then $[L,H]+[H,L]$ is an ideal of $L$;
\item[\upshape(vi)] if $H$ is an ideal of $L$, then $[\gamma_{j}(H),\gamma_{k}(H)]\leq\gamma_{j+k}(H)$ for every pair of positive integers $j,k$;
\item[\upshape(vii)] if $H$ is an ideal of $L$, then $\gamma_{j}(H)$ is an ideal of $L$ for each positive integer $j$; in particular, $\gamma_{j}(L)$ is an ideal of $L$ for each positive integer $j$;
\item[\upshape(viii)] if $H$ is an ideal of $L$, then $\gamma_{j}(\gamma_{k}(H))\leq\gamma_{jk}(H)$ for all positive integers $j,k$.
\end{enumerate}
\end{prop}
\pf
(i) is obvious.

(ii) If $u,v\in[H,H]$, then
$$u=\sum_{1\leq j\leq n}[x_{j},y_{j}]$$
where $x_{j},y_{j}\in H$, $1\leq j\leq n$,
$$v=\sum_{1\leq j\leq t}[a_{j},b_{j}]$$
where $a_{j},b_{j}\in H$, $1\leq j\leq t$. Consider an element $[[x_{j},y_{j}],z]$ where $z\in L$. We have:
$$[[x_{j},y_{j}],z]=[x_{j},[y_{j},z]]-[y_{j},[x_{j},z]].$$
Since $H$ is an ideal of $L$, $[y_{j},z],\ [x_{j},z]\in H$, so that
$$[x_{j},[y_{j},z]],\ [y_{j},[x_{j},z]]\in[H,H].$$
In particular, if $z=v$, then $[u,v]\in [H,H]$, which implies that $[H,H]$ is a subring of $L$. Further,
$$[z,[x_{j},y_{j}]]=[[z,x_{j}],y_{j}]+[x_{j},[z,y_{j}]].$$
Since $H$ is an ideal of $L$, $[z,x_{j}],\ [z,y_{j}]\in H$, so that
$$[[z,x_{j}],y_{j}],\ [x_{j},[z,y_{j}]]\in[H,H].$$
It follows that $[H,H]$ is an ideal of $L$.

(iii) If $u,v\in[L,H]$, then
$$u=\sum_{1\leq j\leq n}[x_{j},y_{j}]$$
where $x_{j}\in L$, $y_{j}\in H$, $1\leq j\leq n$,
$$v=\sum_{1\leq m\leq t}[a_{m},b_{m}]$$
where $a_{m}\in L$, $b_{m}\in H$, $1\leq m\leq t$. We have:
\begin{align*}
[u,v]&=\left[\sum_{1\leq j\leq n}[x_{j},y_{j}],\sum_{1\leq m\leq t}[a_{m},b_{m}]\right]\\
&=\sum_{1\leq j\leq n;\ 1\leq m\leq t}[[x_{j},y_{j}],[a_{m},b_{m}]].
\end{align*}
Since $H$ is an ideal, $[a_{m},b_{m}]\in H$ for each $1\leq m\leq t$. It follows that $[[x_{j},y_{j}],[a_{m},b_{m}]]\in[L,H]$ for all $1\leq j\leq n$, $1\leq m\leq t$. Thus, $[u,v]\in[L,H]$.

(iv) If $u,v\in[H,L]$, then
$$u=\sum_{1\leq j\leq n}[x_{j},y_{j}]$$
where $x_{j}\in H$, $y_{j}\in L$, $1\leq j\leq n$,
$$v=\sum_{1\leq m\leq t}[a_{m},b_{m}]$$
where $a_{m}\in H$, $b_{m}\in L$, $1\leq m\leq t$. We have:
\begin{align*}
[u,v]&=\left[\sum_{1\leq j\leq n}[x_{j},y_{j}],\sum_{1\leq m\leq t}[a_{m},b_{m}]\right]\\
&=\sum_{1\leq j\leq n;\ 1\leq m\leq t}[[x_{j},y_{j}],[a_{m},b_{m}]].
\end{align*}
Since $H$ is an ideal, $[x_{j},y_{j}]\in H$ for each $1\leq j\leq n$. It follows that $[[x_{j},y_{j}],[a_{m},b_{m}]]\in[H,L]$ for all $1\leq j\leq n$, $1\leq m\leq t$. Thus, $[u,v]\in[L,H]$.

(v) Let $h_{1},h_{2},h_{3},h_{4}\in H$, $x_{1},x_{2},x_{3},x_{4}\in L$. Then we have:
\begin{gather*}
[[h_{1},x_{1}]+[x_{2},h_{2}],[h_{3},x_{3}]+[x_{4},h_{4}]]\\
=[[h_{1},x_{1}],[h_{3},x_{3}]]+[[h_{1},x_{1}],[x_{4},h_{4}]]+[[x_{2},h_{2}],[h_{3},x_{3}]]\\
+[[x_{2},h_{2}],[x_{4},h_{4}]]=y.
\end{gather*}
Since $H$ is an ideal, $y\in[H,L]+[L,H]$. It follows that $[L,H]+[H,L]$ is a subring of $L$. If $u\in[L,H]$ (respectively, $v\in[H,L]$), then
$$u=\sum_{1\leq j\leq n}[x_{j},y_{j}]$$
where $x_{j}\in L$, $y_{j}\in H$, $1\leq j\leq n$ (respectively,
$$v=\sum_{1\leq m\leq t}[a_{m},b_{m}]$$
where $a_{m}\in H$, $b_{m}\in L$, $1\leq m\leq t$).

Consider the elements $[[x_{j},y_{j}],z]$ and $[[a_{m},b_{m}],z]$, where $z\in L$. Since $H$ is an ideal, $[x_{j},y_{j}],\ [a_{m},b_{m}]\in H$, so that $[[x_{j},y_{j}],z],\ [[a_{m},b_{m}],z]\in[H,L]$. Similarly, $[z,[x_{j},y_{j}]],\ [z,[a_{m},b_{m}]]\in[L,H]$. It follows that $[u+v,z],\ [z,u+v]\in[L,H]+[H,L]$, which proves that $[L,H]+[H,L]$ is an ideal of $L$.

(vi) We will use the induction on $j$. For $j=1$, the result follows from definition. Now, suppose that $j>1$ and we have already proved the inclusion
$$[\gamma_{m}(H),\gamma_{k}(H)]\leq\gamma_{m+k}(H)$$
for all $m<j$. We have: $\gamma_{j}(H)=[H,\gamma_{j-1}(H)]$. Choose the arbitrary elements $x\in H$, $y\in\gamma_{j-1}(H)$, $z\in\gamma_{k}(H)$. We have
$$[[x,y],z]=[x,[y,z]]-[y,[x,z]].$$
Since $[y,z]\in[\gamma_{j-1}(H),\gamma_{k}(H)]$, the induction hypothesis implies that $[y,z]\in\gamma_{j-1+k}(H)$, so that $[x,[y,z]]\in[H,\gamma_{j-1+k}(H)]=\gamma_{j+k}(H)$. Further,
$$[y,[x,z]]\in[\gamma_{j-1}(H),[H,\gamma_{k}(H)]]=[\gamma_{j-1}(H),\gamma_{k+1}(H)].$$
By the induction hypothesis,
$$[\gamma_{j-1}(H),\gamma_{k+1}(H)]\leq\gamma_{j-1+k+1}(H)=\gamma_{j+k}(H).$$
Hence $[y,[x,z]]\in\gamma_{j+k}(H)$, which proves the inclusion
$$[\gamma_{j}(H),\gamma_{k}(H)]\leq\gamma_{j+k}(H).$$

(vii) Again, we proceed by induction on $j$. For $j=1$, the result follows directly from the definition. Assume now that $j>1$ and that $\gamma_{m}(H)$ is an ideal of $L$ for all $m<j$. We have: $\gamma_{j}(H)=[H,\gamma_{j-1}(H)]$. Choose the arbitrary elements $x\in H$, $y\in\gamma_{j-1}(H)$, $z\in L$. We have
$$[[x,y],z]=[x,[y,z]]-[y,[x,z]].$$
By the induction hypothesis, $[y,z]\in\gamma_{j-1}(H)$, so that
$$[x,[y,z]]\in[H,\gamma_{j-1}(H)]=\gamma_{j}(H).$$
Since $H$ is an ideal of $L$, $[x,z]\in H$, so that
$$[y,[x,z]]\in[\gamma_{j-1}(H),H]=[\gamma_{j-1}(H),\gamma_{1}(H)]\leq\gamma_{j}(H).$$
Hence $[[x,y],z]\in\gamma_{j}(H)$. Further,
$$[z,[x,y]]=[[z,x],y]+[x,[z,y]].$$
Since $H$ is an ideal, $[z,x]\in H$, so $[[z,x],y]\in[H,\gamma_{j-1}(H)]=\gamma_{j}(H)$. By the induction hypothesis, $[z,y]\in\gamma_{j-1}(H)$, so that $[x,[z,y]]\in[H,\gamma_{j-1}(H)]=\gamma_{j}(H)$.

(viii) We will use the induction on $j$. For $j=1$, we have $\gamma_{1}(\gamma_{k}(H))=\gamma_{k}(H)$, so the statement is trivial. Assume now that $j>1$ and that the inclusion $\gamma_{m}(\gamma_{k}(H))\leq\gamma_{mk}(H)$ holds for all $m<j$. Then
$$\gamma_{j}(\gamma_{k}(H))=[\gamma_{k}(H),\gamma_{j-1}(\gamma_{k}(H))]\leq[\gamma_{k}(H),\gamma_{jk-k}(H)].$$
Using (vi), we obtain the inclusion
$$\gamma_{j}(\gamma_{k}(H))\leq\gamma_{k+jk-k}(H)=\gamma_{jk}(H).$$
\qed

By Proposition~\ref{P4}, $[L,L]$ is an ideal of a Leibniz ring $L$; this ideal is called the \textit{derived ideal} of $L$.

We remark that if $A,B$ are ideals of a Leibniz ring $L$, then, in general, $[A,B]$ need not be an ideal. For Leibniz algebras, this was shown by D.~Barnes~\cite{BD2013}, and it is not hard to construct a similar example for Leibniz rings.

Every Leibniz ring has one specific ideal. Denote by $Leib(L)$ the subgroup of the additive group of $L$ generated by the elements $[a,a]$, $a\in L$. From this definition, it follows that the derived ideal $[L,L]$ contains $Leib(L)$.

\begin{prop}\label{P5}
Let $L$ be a Leibniz ring. Then $Leib(L)$ is an ideal of $L$ such that $L/Leib(L)$ is a Lie ring. Moreover, if $H$ is an ideal of $L$ such that $L/H$ is a Lie ring, then $Leib(L)\leq H$.
\end{prop}
\pf
We have
$$[a,[a,x]]=[[a,a],x]+[a,[a,x]],$$
so
$$[[a,a],x]=0.$$
Furthermore,
\begin{align*}
[x+[a,a],x+[a,a]]&=[x,x]+[x,[a,a]]\\
&+[[a,a],x]+[[a,a],[a,a]]\\
&=[x,x]+[x,[a,a]].
\end{align*}
It follows that
$$[x,[a,a]]=[x+[a,a],x+[a,a]]-[x,x]\in Leib(L).$$
Hence $Leib(L)$ is an ideal of $L$.

Put $K=Leib(L)$. Then in factor-ring $L/K$ we have
$$[a+K,a+K]=[a,a]+K=K$$
for each $a\in L$. By Proposition~\ref{P3}, it follows that $L/K$ is a Lie ring.

Now, let $H$ be an ideal of $L$ such that $L/H$ is a Lie ring. Then
$$H=[a+H,a+H]=[a,a]+H$$
for every $a\in L$, which means $[a,a]\in H$. Thus $Leib(L)\leq H$.
\qed

The ideal $Leib(L)$ is called the \textit{Leibniz kernel} of a ring $L$. Note the following important property of $Leib(L)$, which has been proved in Proposition~\ref{P5}.

\begin{prop}\label{P6}
Let $L$ be a Leibniz ring. Then $[Leib(L),L]=\langle0\rangle$. In particular, $Leib(L)$ is an abelian ideal.
\end{prop}

Let $A$ be an additive abelian group and $n$ a positive integer. Define
$$\Lambda_{n}(A)=\{a\in A|\ na=0\}.$$
If $A$ is an additive abelian $p$-group, where $p$ is a prime, then define
$$\Omega_{n}(A)=\{a\in A|\ p^{n}a=0\}.$$

\begin{prop}\label{P7}
Let $L$ be a Leibniz ring. Then $\Lambda_{n}(L)$ is an ideal of $L$ for every positive integer $n$.
\end{prop}
\pf
Let $a\in\Lambda_{n}(A)$ and $x\in L$ be arbitrary. Since $na=0$, we have
$$n[a,x]=[na,x]=0=[x,na]=n[x,a].$$
Thus $[a,x],[x,a]\in\Lambda_{n}(L)$, and therefore $\Lambda_{n}(L)$ is an ideal of $L$.
\qed

\begin{cor}\label{C8}
Let $L$ be a Leibniz ring, $p$ a prime, and $P$ the maximal $p$-subgroup of the additive group of $L$. Then $\Omega_{n}(P)$ is an ideal of $L$ for every positive integer $n$.
\end{cor}

\begin{cor}\label{C9}
Let $L$ be a Leibniz ring and $p$ a prime. Then the maximal $p$-subgroup of the additive group of $L$ is an ideal of $L$.
\end{cor}

\begin{cor}\label{C10}
Let $L$ be a Leibniz ring and $\pi$ a set of primes. Then the maximal $\pi$-subgroup of the additive group of $L$ is an ideal of $L$.
\end{cor}

\begin{cor}\label{C11}
Let $L$ be a Leibniz ring. Then the periodic part of the additive group of $L$ is an ideal of $L$.
\end{cor}

Let $\pi$ be the set of primes and $S$ the maximal $\pi$-subgroup of the additive group of a Leibniz ring $L$. By Corollary~\ref{C10}, $S$ is an ideal of $L$. We will say that $S$ is the \textit{maximal $\pi$-ideal} of $L$.

\begin{cor}\label{C12}
Let $L$ be a Leibniz ring and $T$ the periodic part of the additive group of $L$. Then
$$T=\bigoplus_{p\in\Pi(L)}S_{p},$$
where $S_{p}$ is the maximal $p$-ideal of $L$.
\end{cor}

\begin{prop}\label{P13}
Let $L$ be a Leibniz ring. Then for every positive integer $n$, $nL=\langle na|\ a\in L\rangle$ is an ideal of $L$.
\end{prop}
\pf
Let $a\in nL$ and $x\in L$ be arbitrary. Then $a=nb$ for some $b\in L$. We have
\begin{align*}
[a,x]&=[nb,x]=n[b,x]\in nL,\\
[x,a]&=[x,nb]=n[x,b]\in nL.
\end{align*}
Thus $[a,x],[x,a]\in nL$, and therefore $nL$ is an ideal of $L$.
\qed

Let $L$ be a Leibniz ring and put
$$\alpha(L)=\{a\in L|\ [x,a]=-[a,x]\ \mbox{for all }x\in L\}.$$
The subset $\alpha(L)$ is called the \textit{anticenter} of the ring $L$.

\begin{prop}\label{P14}
Let $L$ be a Leibniz ring. Then the anticenter $\alpha(L)$ of $L$ is an ideal of $L$.
\end{prop}
\pf
Let $a,b\in\alpha(L)$ and $x\in L$ be arbitrary. Then
$$[a-b,x]=[a,x]-[b,x]=-[x,a]+[x,b]=-([x,a]-[x,b])=-[x,a-b].$$
Thus $a-b\in\alpha(L)$, and hence $\alpha(L)$ is a subgroup of the additive group of $L$. Now let $y\in L$. Then
\begin{align*}
[[a,y],x]&=[a,[y,x]]-[y,[a,x]]\\
&=-[[y,x],a]-[y,[a,x]]\\
&=-[y,[x,a]]+[x,[y,a]]-[y,[a,x]]\\
&=-[y,[x,a]]-[x,[a,y]]+[y,[x,a]]\\
&=-[x,[a,y]].
\end{align*}
It follows that $[a,y]\in\alpha(L)$. Since $[y,a]=-[a,y]\in\alpha(L)$, we conclude that $\alpha(L)$ is an ideal of $L$.
\qed

Proposition~\ref{P1} implies the following

\begin{cor}\label{C15}
Let $L$ be a symmetric Leibniz ring. Then the anticenter of $L$ contains $[L,L]$. In particular, if $L=[L,L]$, then $L$ is a Lie ring.
\end{cor}

Let $L$ be a Leibniz ring and put
\begin{align*}
\zeta^{\mathrm{left}}(L)&=\{x\in L|\ [x,y]=0\ \mbox{for all}\ y\in L\},\\
\zeta^{\mathrm{right}}(L)&=\{x\in L|\ [y,x]=0\ \mbox{for all}\ y\in L\},\\
\zeta(L)&=\{x\in L|\ [x,y]=[y,x]=0\ \mbox{for all}\ y\in L\}.
\end{align*}
The subset $\zeta^{\mathrm{left}}(L)$ is called the \textit{left center} of $L$, the subset $\zeta^{\mathrm{right}}(L)$ is called the \textit{right center} of $L$, and the subset $\zeta(L)$ is called the \textit{center} of $L$.

\begin{prop}\label{P16}
Let $L$ be a Leibniz ring. Then $\zeta^{\mathrm{right}}(L)$ is a subring of $L$, while $\zeta^{\mathrm{left}}(L)$ and $\zeta(L)$ are ideals of $L$.
\end{prop}
\pf
The fact that $\zeta^{\mathrm{left}}(L)$ and $\zeta^{\mathrm{right}}(L)$ are subrings of $L$ is almost obvious.

Let $x,y\in L$ and $z\in\zeta^{\mathrm{left}}(L)$. Then
$$[[z,x],y]=[0,y]=0$$
and
$$[[x,z],y]=[x,[z,y]]-[z,[x,y]]=[x,0]-0=0.$$
These equalities show that $\zeta^{\mathrm{left}}(L)$ is an ideal of $L$. The fact that $\zeta(L)$ is an ideal of $L$ is clear.
\qed

We note that, in general, the left and right centers need not coincide. Concrete examples can be found in~\cite[Chapter~1]{KPS2024}.

\begin{prop}\label{P17}
Let $L$ be a symmetric Leibniz ring. Then $\zeta^{\mathrm{right}}(L)$, $\zeta^{\mathrm{left}}(L)$, and $\zeta(L)$ are ideals of $L$, and moreover
\begin{align*}
\zeta^{\mathrm{right}}(L)/\zeta(L)&\leq\zeta(L/\zeta(L)),\\
\zeta^{\mathrm{left}}(L)/\zeta(L)&\leq\zeta(L/\zeta(L)).
\end{align*}
\end{prop}
\pf
Let $a\in\zeta^{\mathrm{right}}(L)$, $x,y\in L$. Then
$$[x,a]=0\in\zeta^{\mathrm{right}}(L).$$
Furthermore,
$$[y,[a,x]]=[[y,a],x]+[a,[y,x]]=[a,[y,x]].$$
By Proposition~\ref{P1} we have
$$[a,[y,x]]=-[[y,x],a]=0,$$
and hence $[a,x]\in\zeta^{\mathrm{right}}(L)$. Thus $\zeta^{\mathrm{right}}(L)$ is an ideal of $L$.

Let $b=[a,x]$. From the above we know that $b\in\zeta^{\mathrm{right}}(L)$, so $[y,b]=0$ for every $y\in L$. Since $b\in[L,L]$, Corollary~\ref{C15} implies that
$$[b,y]=-[y,b]=0.$$
Therefore $b\in\zeta(L)$. In other words,
$$[\zeta^{\mathrm{right}}(L),L]\leq\zeta(L).$$
On the other hand,
$$[L,\zeta^{\mathrm{right}}(L)]=\langle0\rangle\leq\zeta(L)$$
which implies that
$$\zeta^{\mathrm{right}}(L)/\zeta(L)\leq\zeta(L/\zeta(L)).$$
By similar arguments we obtain that
$$\zeta^{\mathrm{left}}(L)/\zeta(L)\leq\zeta(L/\zeta(L)).$$
\qed

Consider now the structure of Leibniz rings whose additive groups are ``small''. As will be seen, these results already demonstrate how much more complicated and diverse the situation with Leibniz rings is.

\section{Structure of Leibniz rings with cyclic additive group}
Let $L$ be a Leibniz ring. The order of an element $x$, regarded as an element of the additive group of $L$, is called the \textit{additive order} of $x$.

Suppose that $L$ is a Leibniz ring whose additive group is cyclic, and let $a$ be a generator of the additive group of $L$. Then
$$L=\{ka|\ k\in\mathbb{Z}\}.$$

If $L$ is a Lie ring, then the equality $[a,a]=0$ implies that
$$[ka,ta]=kt[a,a]=0.$$
Thus, in this case, $L$ is a ring with zero multiplication.

Therefore, suppose that $L$ is not a Lie ring. Then $Leib(L)=K$ is a non-zero ideal of $L$. Since $K$ is an abelian ideal by Proposition~\ref{P6}, we have $K\neq L$. The additive group of $K$ is generated by an element $kg$ for some positive integer $k$. The fact that $kg\in K$, together with Proposition~\ref{P6}, implies that $[kg,g]=0$. It follows that
$$0=[kg,g]=k[g,g].$$
In other words, the element $[g,g]$ has finite order in the additive group of $L$. Thus we again obtain that $[g,g]=0$.

Suppose now that the additive group of $L$ is a finite cyclic group. By Corollary~\ref{C12}, we may assume that the additive group of $L$ is a $p$-group for some prime $p$. Let $|L|=p^{m}$ for some positive integer $m$. Since $\langle0\rangle\neq K\neq L$, we have $p^{t}g\in K$ for some $0<t<m$. Then
$$0=[p^{t}g,g]=p^{t}[g,g].$$
It follows that the element $[g,g]$ has additive order $p^{s}$ for some $s\leq t$.

Let $x,y\in L$. Then $x=kg$, $y=jg$ for some non-negative integers $k,j$. We have
$$[x,y]=[kg,jg]=kj[g,g].$$
Hence $[L,L]\leq\langle[g,g]\rangle$. On the other hand, since $[g,g]\in[L,L]$, we obtain $\langle[g,g]\rangle\leq[L,L]$. Thus
$$[L,L]=\langle[g,g]\rangle.$$

The additive factor-group $\langle g\rangle/\langle[g,g]\rangle$ has order
$$\frac{p^{m}}{p^{s}}=p^{m-s}.$$
Therefore $[g,g]=p^{m-s}g$ and
$$0=[[g,g],g]=[p^{m-s}g,g]=p^{m-s}[g,g]=[g,p^{m-s}g]=[g,[g,g]].$$
It follows that
$$[[g,g],kg]=k[[g,g],g]=0$$
and
$$[kg,[g,g]]=k[g,[g,g]]=0.$$
These equalities show that $[g,g]\in\zeta(L)$ and hence $\langle[g,g]\rangle\leq\zeta(L)$. Moreover, since $p^{m-s}[g,g]=0$, it follows that $p^{s}$ divides $p^{m-s}$. Hence $s\leq m-s$, i.e. $2s\leq m$.

It is not hard to see that $x\in\zeta^{\mathrm{left}}(L)$ if and only if $[x,g]=0$. Let $x=jg$, where $j=p^{t}r$ and $p$ does not divide $r$. Then
$$0=[x,g]=[jg,g]=j[g,g]=jp^{m-s}g=p^{t}rp^{m-s}g=rp^{m-s+t}g.$$
It follows that $m-s+t\geq m$, hence $t\geq s$. Thus we obtain $\zeta^{\mathrm{left}}(L)=\langle p^{s}g\rangle$. By the same arguments, $\zeta^{\mathrm{right}}(L)=\langle p^{s}g\rangle$ and hence $\zeta(L)=\langle p^{s}g\rangle$.

Conversely, suppose that $L=\langle g\rangle$ is an additive cyclic group of order $p^{m}$, where $m$ is a positive integer. Put $c=p^{m-s}g$, $2s\leq m$. Then the element $c$ has order $p^{s}$. Define a multiplication $[,]$ on $L$ by the rule: if $x=n_{1}g$, $y=n_{2}g$ with $n_{1},n_{2}$ are non-negative integers, then
$$[x,y]=n_{1}n_{2}c.$$
In particular, $[g,g]=c$. Moreover,
$$[g,c]=[g,p^{m-s}g]=p^{m-s}[g,g]=p^{m-s}c=0,$$
since $m-s\geq s$. Similarly,
$$0=p^{m-s}c=p^{m-s}[g,g]=[p^{m-s}g,g]=[c,g].$$
Thus $c\in\zeta(L)$. Consequently, $[x,y]\in\zeta(L)$ for all $x,y\in L$. If $x,y,z\in L$ are arbitrary, then by the above observations
$$[x,y],[x,z],[y,z]\in\zeta(L),$$
so that
$$0=[x,[y,z]]=[[x,y],z]+[y,[x,z]]=0.$$
It follows that $L$ is a Leibniz ring. Thus we obtain the following type of Leibniz ring:
\begin{gather*}
L_{1}=\langle g\rangle,\\
|g|=p^{m},\ [g,g]=c,\ |c|=p^{s},\ \mbox{where}\ 2s\leq m,\ [g,c]=[c,g]=0.
\end{gather*}
Here
$$Leib(L_{1})=[L_{1},L_{1}]=\langle p^{m-s}g\rangle,\ \zeta^{\mathrm{right}}(L_{1})=\zeta^{\mathrm{left}}(L_{1})=\zeta(L_{1})=\langle p^{s}g\rangle.$$

\section{Structure of Leibniz rings whose additive group is elementary abelian of order $p^{2}$}
Let $L$ be a Leibniz ring and suppose that the additive group of $L$ is an elementary abelian $p$-group of order $p^{2}$. Assume that $L$ is not a Lie ring. As above, we then have $Leib(L)=K$ is a non-zero ideal of $L$, and $K\neq L$. It follows that $K$ has prime order $p$. Since $L/K$ is a Lie ring generated by one element, we obtain that $L/K$ is a ring with zero multiplication. Because $L$ is not a Lie ring, there exists an element $a\in L$ such that $[a,a]\neq0$. Put
$$[a,a]=b.$$
Then $b\in K$, hence $K=\langle b\rangle$. With this choice we have $[b,a]=0$. Since $K$ is an ideal, there exists $0\leq k<p$ such that $[a,c]=kc$. Suppose $k\neq 0$, and let $t$ be a positive integer such that
$$kt\equiv1(\mathrm{mod}\ p).$$
Put $c=ta$. Then
$$[c,b]=[ta,b]=t[a,b]=tkb=b.$$
Also,
$$[c,c]=t^{2}[a,a]=t^{2}b=d,$$
and
$$[c,d]=[c,t^{2}b]=t^{2}[c,b]=t^{2}b=d.$$
Thus, with this choice we have
$$L=\langle c\rangle\oplus\langle d\rangle.$$
Hence we obtain the following two types of Leibniz rings:
\begin{gather*}
L_{2}=\langle a\rangle\oplus\langle b\rangle,\\
pa=pb=0,\ [a,a]=b,\ [b,a]=[a,b]=[b,b]=0.
\end{gather*}
Here
$$Leib(L_{2})=[L_{2},L_{2}]=\zeta(L_{2})=\langle b\rangle.$$
The second type is the following:
\begin{gather*}
L_{3}=\langle c\rangle\oplus\langle d\rangle,\\
[c,c]=[c,d]=d,\ pc=pd=0,\ [d,c]=[d,d]=0.
\end{gather*}
Here
$$Leib(L_{3})=[L_{3},L_{3}]=\zeta^{\mathrm{left}}(L_{3})=\langle d\rangle,\ \zeta^{\mathrm{right}}(L_{3})=\zeta(L_{3})=\langle0\rangle.$$

\section{Structure of Leibniz rings whose additive group is a direct sum of a cyclic subgroup of order $p^{2}$ and a subgroup of order $p$}
Let $L$ be a Leibniz ring and suppose that the additive group of $L$ has order $p^{3}$ and is a direct sum of two cyclic subgroups:
$$L=\langle a\rangle\oplus\langle b\rangle,$$
where the element $a$ has additive order $p^{2}$, and $b$ has additive order $p$. As before, we assume that $L$ is not a Lie ring. Then, as above, $Leib(L)=K$ is a non-zero ideal of $L$ and $K\neq L$.

Suppose first that $K$ has prime order $p$. Then
$$K=\langle u\rangle,\ u=t_{1}a+t_{2}b,\ 0\leq t_{1},t_{2}<p.$$
If $t_{2}\neq0$, then $\langle a\rangle\cap\langle u\rangle=\langle0\rangle$ and hence
$$L=\langle a\rangle\oplus\langle u\rangle.$$
Therefore, without loss of generality, we may assume that $u=b$.

It is not hard to prove that a Lie ring whose additive group is cyclic has zero multiplication. It follows that
$$[a,a],[a,b]\in\langle b\rangle.$$
Without loss of generality, we may assume that $[a,a]=b$ and $[a,b]=\alpha b$ where $0\leq\alpha<p$. If $\alpha=0$, then we obtain the following Leibniz ring:
\begin{gather*}
L_{4}=\langle a\rangle\oplus\langle b\rangle,\\
p^{2}a=pb=0,\ [a,a]=b,\ [b,a]=[a,b]=[b,b]=0.
\end{gather*}
Here
$$Leib(L_{4})=[L_{4},L_{4}]=\zeta^{\mathrm{left}}(L_{4})=\zeta^{\mathrm{right}}(L_{4})=\zeta(L_{4})=\langle b\rangle.$$
Conversely, since $L_{4}/\zeta(L_{4})$ is abelian, it is not hard to see that this ring is indeed a Leibniz ring.

Suppose now that $\alpha\neq0$. Let $\beta$ be a positive integer such that $\alpha\beta\equiv1(\mathrm{mod}\ p)$, and put $c=\beta a$. Then
$$[c,b]=[\beta a,b]=\beta[a,b]=\beta\alpha b=b.$$
We also have
$$[c,c]=\beta^{2}[a,a]=\beta^{2}b=d,$$
and
$$[c,d]=[c,\beta^{2}b]=\beta^{2}[c,b]=\beta^{2}b=d.$$
Thus we obtain the following type of Leibniz ring:
\begin{gather*}
L_{5}=\langle c\rangle\oplus\langle d\rangle,\\
p^{2}c=pd=0,\ [c,c]=[c,d]=d,\ [d,c]=[d,d]=0.
\end{gather*}
Here
$$Leib(L_{5})=[L_{5},L_{5}]=\zeta^{\mathrm{left}}(L_{5})=\langle d\rangle,\ \zeta^{\mathrm{right}}(L_{5})=\zeta(L_{5})=\langle0\rangle.$$

Conversely, we show that such a ring really is a Leibniz ring. Let $x,y,z$ be arbitrary elements of $L_{5}$,
\begin{align*}
x&=\lambda_{1}c+\mu_{1}d,\\
y&=\lambda_{2}c+\mu_{2}d,\\
z&=\lambda_{3}c+\mu_{3}d,
\end{align*}
$0\leq\lambda_{1},\lambda_{2},\lambda_{3}<p^{2}$ and $0\leq\mu_{1},\mu_{2},\mu_{3}<p$. We have
\begin{align*}
[x,y]&=[\lambda_{1}c+\mu_{1}d,\lambda_{2}c+\mu_{2}d]\\
&=\lambda_{1}\lambda_{2}[c,c]+\mu_{1}\lambda_{2}[d,c]+\lambda_{1}\mu_{2}[c,d]+\mu_{1}\mu_{2}[d,d]\\
&=\lambda_{1}\lambda_{2}d+\lambda_{1}\mu_{2}d=(\lambda_{1}\lambda_{2}+\lambda_{1}\mu_{2})d,\\
[x,z]&=(\lambda_{1}\lambda_{3}+\lambda_{1}\mu_{3})d,\\
[y,z]&=(\lambda_{2}\lambda_{3}+\lambda_{2}\mu_{3})d.
\end{align*}
Therefore,
\begin{align*}
[x,[y,z]]&=[\lambda_{1}c+\mu_{1}d,(\lambda_{2}\lambda_{3}+\lambda_{2}\mu_{3})d]\\
&=\lambda_{1}(\lambda_{2}\lambda_{3}+\lambda_{2}\mu_{3})[c,d]\\
&=\lambda_{1}\lambda_{2}(\lambda_{3}+\mu_{3})d,\\
[[x,y],z]&=[(\lambda_{1}\lambda_{2}+\lambda_{1}\mu_{2})d,\lambda_{3}c+\mu_{3}d]\\
&=0,\\
[y,[x,z]]&=[\lambda_{2}c+\mu_{2}d,(\lambda_{1}\lambda_{3}+\lambda_{1}\mu_{3})d]\\
&=\lambda_{2}(\lambda_{1}\lambda_{3}+\lambda_{1}\mu_{3})[c,d]\\
&=\lambda_{2}\lambda_{1}(\lambda_{3}+\mu_{3})d.
\end{align*}
Hence
$$[x,[y,z]]=[[x,y],z]+[y,[x,z]],$$
which shows that $L_{5}$ indeed is a Leibniz ring.

Suppose now that $K$ has order $p^{2}$. Then either $K$ is cyclic or $K$ is elementary abelian. Assume that $K$ is cyclic. Without loss of generality, we may suppose that $K=\langle a\rangle$. In this case we have
$$[b,b],[b,a]\in K.$$
Since
$$0=[0,a]=[pb,a]=p[b,a],$$
it follows that $p[b,a]=0$. Similarly, $0=p[b,b]$. It is not hard to show that in this case $[L,L]\leq pK$. Hence $|[L,L]|\leq p$. On the other hand, since $K\leq[L,L]$, we obtain a contradiction. This contradiction shows that additive group of $K$ must be elementary abelian.

Proposition~\ref{P13} implies that $pL$ is an ideal of $L$. We note that $pL=\langle pa\rangle$. Then $|pL|=p$, in particular, $pL\neq K$. It follows that $L/pL$ is not a Lie ring. In particular, $Leib(L/pL)\neq L/pL$, so that $|Leib(L/pL)|=p$. The inclusion $K/pL\leq Leib(L/pL)$ implies that
$$K/pL=Leib(L/pL).$$
By the argument above, $K$ does not contain $\langle a\rangle$. Therefore, without loss of generality, we may assume that
$$K/pL=(\langle b\rangle+pL)/pL,$$
so that
$$K=\langle pa\rangle\oplus\langle b\rangle.$$
Put $c=pa$. We have
$$0=p[a,a]=[a,pa].$$
Since $[b,c]=[c,b]=0$, it follows that $c\in\zeta(L)$.

Suppose that $[a,b]=0$. Let $x$ be an arbitrary element of $L$,
$$x=\lambda a+\mu b,$$
$0\leq\lambda<p^{2}$, $0\leq\mu<p$. Then
\begin{align*}
[x,x]&=[\lambda a+\mu b,\lambda a+\mu b]\\
&=\lambda^{2}[a,a]+\lambda\mu[a,b]+\mu\lambda[b,a]+\mu^{2}[b,b]\\
&=\lambda^{2}[a,a]+\lambda\mu[a,b]=\lambda^{2}[a,a].
\end{align*}
Thus, in this case $K=\langle[a,a]\rangle$ is cyclic. This contradicts the assumption that $K$ is elementary abelian. Therefore, we conclude that $[a,b]\neq0$.

Let $[a,b]=\beta b+\gamma c$, where $0<\beta<p$, $0\leq\gamma<p$. Without loss of generality we may suppose that $\beta=1$ (otherwise we replace $b$ with $\beta^{-1}b$). Put $a_{2}=b+\gamma c$. Then
$$[a,a_{2}]=[a,b+\gamma c]=[a,b]+\gamma[a,c]=[a,b]=b+\gamma c=a_{2}.$$
The equalities
$$[a_{2},a]=[a_{2},b]=[b,a_{2}]=[a_{2},a_{2}]=0$$
shows that $\langle a_{2}\rangle$ is an ideal of $L$. The additive group of the factor-ring $L/\langle a_{2}\rangle$ is cyclic of order $p^{2}$. On the other hand, since $\langle a_{2}\rangle\neq K$, the factor-ring $L/\langle a_{2}\rangle$ is not a Lie ring. Taking into account the structure of a Leibniz ring with cyclic additive group (described above), we obtain $[a,a]=c+\beta_{1}a_{2}$. Now put $a_{1}=a-\beta_{1}a_{2}$. Then
$$[a_{1},a_{1}]=[a-\beta_{1}a_{2},a-\beta_{1}a_{2}]=[a,a]-\beta_{1}[a,a_{2}]=c+\beta_{1}a_{2}-\beta_{1}a_{2}=c.$$
Note also that
$$pa_{1}=p(a-\beta_{1}a_{2})=pa=c.$$
Thus we come to the following type of Leibniz ring:
\begin{gather*}
L_{6}=\langle a_{1}\rangle\oplus\langle a_{2}\rangle,\\
p^{2}a_{1}=pa_{2}=0,\ [a_{1},a_{1}]=pa_{1},\ [a_{1},a_{2}]=a_{2},\ [a_{2},a_{1}]=[a_{2},a_{2}]=0,\\
[pa_{1},a_{1}]=[a_{1},pa_{1}]=[pa_{1},a_{2}]=[a_{2},pa_{1}]=[pa_{1},pa_{1}]=0.
\end{gather*}

\section{Structure of Leibniz rings whose additive group is the direct sum of two infinite cyclic subgroups}
Let $L$ be a Leibniz ring and suppose that the additive group of $L$ is a direct sum of two infinite cyclic subgroups. Again, we assume that $L$ is not a Lie ring. Then, as above, $Leib(L)=K$ is a non-zero ideal of $L$ and $K\neq L$.

Denote by $T/K$ the periodic part of the additive group of $L/K$. By Corollary~\ref{C11}, $T$ is an ideal of $L$. If we suppose that $T=L$, then we obtain $L=\zeta^{\mathrm{right}}(L)$. In particular, this means that $L$ is abelian, which gives a contradiction. Therefore, this contradiction shows that $T\neq L$.

Suppose that $T\neq K$ and let $u\in T$ with $u\not\in K$. Then there exists a positive integer $n$ such that $nu\in K$. It follows that $[nu,x]=0$ for each $x\in L$. Hence $0=[nu,x]=n[u,x]$. Since the additive group of $L$ is torsion-free, we obtain $[u,x]=0$. This holds for every $x\in L$, and therefore $u\in\zeta^{\mathrm{left}}(L)$. Thus, $T\leq\zeta^{\mathrm{left}}(L)$. If we assume that $T\neq\zeta^{\mathrm{left}}(L)$, then $\zeta^{\mathrm{left}}(L)$ has finite index in $L$. Using the above arguments, this leads to a contradiction. Hence, we conclude that $T=\zeta^{\mathrm{left}}(L)$.

Since $r_{0}(L)=2$ and $T\neq\langle0\rangle$, the additive group of $L/T$ is infinite cyclic. In this case we obtain $L=T\oplus\langle b\rangle$ for some element $b$. As shown above, the factor-ring $L/T$ is abelian. Hence $[L,L]\leq T$.

Finally, the equality $r_{0}(L)=2$ implies that the additive group of $T$ is infinite cyclic, that is, $T=\langle a\rangle$ for some $a\in L$.

If $T=K$, then as above we obtain that $K=\langle a\rangle$ and the additive group of $L/K$ is infinite cyclic. Hence,
$$L=\langle a\rangle\oplus\langle b\rangle.$$
From the previous arguments it follows that in this case $K=\zeta^{\mathrm{left}}(L)$.

We have $[b,b]=\beta a$, $[b,a]=\alpha a$ for some positive integers $\alpha,\beta$. The equality $\langle a\rangle=\zeta^{\mathrm{left}}(L)$ implies that $[a,a]=[a,b]=0$. Let $x,y,z\in L$ be arbitrary elements,
\begin{align*}
x&=\lambda_{1}a+\mu_{1}b,\\
y&=\lambda_{2}a+\mu_{2}b,\\
z&=\lambda_{3}a+\mu_{3}b.
\end{align*}
Then
\begin{align*}
[x,y]&=[\lambda_{1}a+\mu_{1}b,\lambda_{2}a+\mu_{2}b]\\
&=\lambda_{1}\lambda_{2}[a,a]+\mu_{1}\lambda_{2}[b,a]+\lambda_{1}\mu_{2}[a,b]+\mu_{1}\mu_{2}[b,b]\\
&=\mu_{1}\lambda_{2}[b,a]+\mu_{1}\mu_{2}[b,b]=\mu_{1}\lambda_{2}\alpha a+\mu_{1}\mu_{2}\beta a\\
&=(\mu_{1}\lambda_{2}\alpha+\mu_{1}\mu_{2}\beta)a,\\
[x,z]&=(\mu_{1}\lambda_{3}\alpha+\mu_{1}\mu_{3}\beta)a,\\
[y,z]&=(\mu_{2}\lambda_{3}\alpha+\mu_{2}\mu_{3}\beta)a.
\end{align*}
Hence,
\begin{align*}
[x,[y,z]]&=[\lambda_{1}a+\mu_{1}b,(\mu_{2}\lambda_{3}\alpha+\mu_{2}\mu_{3}\beta)a]\\
&=[\mu_{1}b,(\mu_{2}\lambda_{3}\alpha+\mu_{2}\mu_{3}\beta)a]\\
&=\mu_{1}(\mu_{2}\lambda_{3}\alpha+\mu_{2}\mu_{3}\beta)[b,a]\\
&=\alpha\mu_{1}(\mu_{2}\lambda_{3}\alpha+\mu_{2}\mu_{3}\beta)a,\\
[[x,y],z]&=[(\mu_{1}\lambda_{2}\alpha+\mu_{1}\mu_{2}\beta)a,\lambda_{3}a+\mu_{3}b]\\
&=0,\\
[y,[x,z]]&=[\lambda_{2}a+\mu_{2}b,(\mu_{1}\lambda_{3}\alpha+\mu_{1}\mu_{3}\beta)a]\\
&=[\mu_{2}b,(\mu_{1}\lambda_{3}\alpha+\mu_{1}\mu_{3}\beta)a]\\
&=\mu_{2}(\mu_{1}\lambda_{3}\alpha+\mu_{1}\mu_{3}\beta)[b,a]\\
&=\alpha\mu_{2}(\mu_{1}\lambda_{3}\alpha+\mu_{1}\mu_{3}\beta)a.
\end{align*}
Hence
$$[x,[y,z]]=[[x,y],z]+[y,[x,z]].$$
Thus we arrive at the following type of Leibniz ring:
\begin{gather*}
L_{7}=\langle a_{1}\rangle\oplus\langle a_{2}\rangle,\\
\mbox{where the additive orders of}\ a_{1}\ \mbox{and}\ a_{2}\ \mbox{are infinite},\\
[a_{1},a_{1}]=[a_{1},a_{2}]=0,\ [a_{2},a_{1}]=\alpha a_{1},\ [a_{2},a_{2}]=\beta a_{1}.
\end{gather*}
Here
\begin{gather*}
\zeta^{\mathrm{left}}(L_{7})=\langle a_{1}\rangle,\ \zeta^{\mathrm{right}}(L_{7})=\zeta(L_{7})=\langle 0\rangle,\\ Leib(L_{7})=\langle\beta a_{1}\rangle,\ [L_{7},L_{7}]=\langle\alpha a_{1}\rangle+\langle\beta a_{1}\rangle.
\end{gather*}

\section{Structure of Leibniz rings whose additive group is the direct sum of infinite and finite cyclic subgroups}
Let $L$ be a Leibniz ring, and suppose that the additive group of $L$ is a direct sum of one infinite cyclic subgroup and one finite cyclic subgroup. Again, we assume that $L$ is not a Lie ring. As above, $Leib(L)=K$ is a non-zero ideal of $L$, and $K\neq L$.

Let $T$ be the periodic part of the additive group of $L$. By Corollary~\ref{C11}, $T$ is an ideal of $L$. The additive group of the factor-ring $L/T$ is infinite cyclic. By the arguments proved above, this factor-ring is abelian. It follows that $[L,L]\leq T$ and $K\leq T$. Put $T=\langle a\rangle$. Since the additive group of $L/T$ is infinite cyclic, $L=T\oplus\langle b\rangle$ for some element $b$ of infinite additive order. By our assumptions, $T=\langle a\rangle$ for some element $a$. Let $|a|=k$. The inclusion $K\leq\langle a\rangle$ implies that $[b,b]=\beta a$ for some integer $\beta$, $0\leq\beta<k$. Let $\mathrm{GCD}(\beta,k)=d$ and $\beta=\beta_{1}\beta_{2}$ where $\mathrm{GCD}(\beta_{2},k)=1$. Therefore, without loss of generality, we may assume that $\beta$ is a divisor of $k$.

We have $[a,a]=\sigma a,\ [b,a]=\alpha_{1}a,\ [a,b]=\alpha_{2}a,\ [b,b]=\beta a$, where $\sigma,\alpha_{1},\alpha_{2},\beta$ are integers satisfying $0\leq\sigma,\alpha_{1},\alpha_{2},\beta\leq k$.

Suppose first that $\sigma=\alpha_{1}=\alpha_{2}=0$, that is $\langle a\rangle\leq\zeta(L)$. Note that $\beta\neq k$, otherwise $L$ would be abelian, contradicting our assumptions. Let $x=\lambda_{1}a+\mu_{1}b$ be an arbitrary element of $L$. We have
\begin{align*}
[x,b]&=[\lambda_{1}a+\mu_{1}b,b]=\lambda_{1}[a,b]+\mu_{1}[b,b]=\mu_{1}\beta a,\\ 
[b,x]&=[b,\lambda_{1}a+\mu_{1}b]=\lambda_{1}[b,a]+\mu_{1}[b,b]=\mu_{1}\beta a.
\end{align*}
We note that $x\in\zeta(L)$ if and only if $[x,b]=[b,x]=0$. It follows that $\mu_{1}\beta a=0$, hence $\mu_{1}=k/\beta$. Therefore, in this case we obtain
$$\zeta(L)=\{\lambda a+\mu(k/\beta)b|\ \lambda,\mu\ \mbox{are integers}\}.$$
Thus we arrive at the following type of Leibniz ring:
\begin{gather*}
L_{8}=\langle a_{1}\rangle\oplus\langle a_{2}\rangle,\\
\mbox{the additive order of}\ a_{1}\ \mbox{is finite},\ |a_{1}|=k,\\
\mbox{the additive order of}\ a_{2}\ \mbox{is infinite},\\
[a_{1},a_{1}]=[a_{1},a_{2}]=[a_{2},a_{1}]=0,[a_{2},a_{2}]=\beta a_{1},\\
\mbox{where}\ \beta\ \mbox{is a divisor of}\ k.
\end{gather*}
Here
\begin{gather*}
\langle a_{1}\rangle=\zeta(L_{8})=\langle a_{1}\rangle\oplus\langle(k/\beta)a_{2}\rangle=\zeta^{\mathrm{left}}(L_{8})=\zeta^{\mathrm{right}}(L_{8}),\\
Leib(L_{8})=\langle\beta a_{1}\rangle=[L_{8},L_{8}].
\end{gather*}
Suppose now that the center of $L$ does not contain $\langle a\rangle$. By the arguments above, in this case we have $[a,[a,a]]=0=[[a,a],a]$. We have
$$[a,[a,a]]=[a,\sigma a]=\sigma[a,a]=\sigma^{2}a,$$
so that $\sigma^{2}\equiv0(\mathrm{mod}\ k)$. Furthermore, we have
\begin{align*}
0&=[[a,a],b]=[\sigma a,b]=\sigma[a,b]=\sigma\alpha_{2}a,\ \mbox{so that}\ \alpha_{2}\sigma\equiv0(\mathrm{mod}\ k),\\
0&=[[b,b],b]=[\beta a,b]=\beta[a,b]=\beta\alpha_{2}a,\ \mbox{so that}\ \alpha_{2}\beta\equiv0(\mathrm{mod}\ k),\\
0&=[[b,b],a]=[\beta a,a]=\beta[a,a]=\beta\sigma a,\ \mbox{so that}\ \beta\sigma\equiv0(\mathrm{mod}\ k),
\end{align*}
and
\begin{align*}
[a,[b,b]]&=[[a,b],b]+[b,[a,b]]\\
&=[\alpha_{2}a,b]+[b,\alpha_{2}a]=\alpha_{2}[a,b]+\alpha_{2}[b,a]\\
&=\alpha_{2}^{2}a+\alpha_{2}\alpha_{1}a,\ \mbox{so that}\ \alpha_{2}^{2}+\alpha_{2}\alpha_{1}\equiv0(\mathrm{mod}\ k)\\
[a,[b,b]]&=[a,\beta a]=\beta[a,a]=\beta\sigma a,\\
[b,[a,a]]&=[[b,a],a]+[a,[b,a]]=[\alpha_{1}a,a]+[a,\alpha_{1}a]=2\sigma\alpha_{1}a,\\
[b,[a,a]]&=[b,\sigma a]=\sigma[b,a]=\sigma\alpha_{1}a,\ \mbox{so that}\ \alpha_{1}\sigma\equiv0(\mathrm{mod}\ k).
\end{align*}

Conversely, suppose that $L=\langle a\rangle\oplus\langle b\rangle$ where $|a|=k$ and $b$ has infinite additive order. Define the multiplication $[,]$ on $L$ by the following rules:
$$[a,a]=\sigma a,\ [b,a]=\alpha_{1}a,\ [a,b]=\alpha_{2}a,\ [b,b]=\beta a$$
where $\sigma,\alpha_{1},\alpha_{2},\beta$ are integers with $0\leq\sigma,\alpha_{1},\alpha_{2},\beta\leq k$. Assume that the following congruences are satisfied:
\begin{align*}
\sigma^{2}&\equiv0(\mathrm{mod}\ k),\ \alpha_{2}\sigma\equiv0(\mathrm{mod}\ k),\alpha_{2}\beta\equiv0(\mathrm{mod}\ k),\\ \beta\sigma&\equiv0(\mathrm{mod}\ k),\alpha_{2}^{2}+\alpha_{2}\alpha_{1}\equiv0(\mathrm{mod}\ k),\ \alpha_{1}\sigma\equiv0(\mathrm{mod}\ k).
\end{align*}
Let $x,y,z\in L$ be arbitrary elements,
$$x=\lambda_{1}a+\mu_{1}b,\ y=\lambda_{2}a+\mu_{2}b,\ z=\lambda_{3}a+\mu_{3}b.$$
Then
\begin{align*}
[x,y]&=[\lambda_{1}a+\mu_{1}b,\lambda_{2}a+\mu_{2}b]\\
&=\lambda_{1}\lambda_{2}[a,a]+\mu_{1}\lambda_{2}[b,a]+\lambda_{1}\mu_{2}[a,b]+\mu_{1}\mu_{2}[b,b]\\
&=\lambda_{1}\lambda_{2}\sigma a+\mu_{1}\lambda_{2}\alpha_{1}a+\lambda_{1}\mu_{2}\alpha_{2}a+\mu_{1}\mu_{2}\beta a\\
&=(\lambda_{1}\lambda_{2}\sigma+\mu_{1}\lambda_{2}\alpha_{1}+\lambda_{1}\mu_{2}\alpha_{2}+\mu_{1}\mu_{2}\beta)a,\\
[x,z]&=(\lambda_{1}\lambda_{3}\sigma+\mu_{1}\lambda_{3}\alpha_{1}+\lambda_{1}\mu_{3}\alpha_{2}+\mu_{1}\mu_{3}\beta)a,\\ 
[y,z]&=(\lambda_{2}\lambda_{3}\sigma+\mu_{2}\lambda_{3}\alpha_{1}+\lambda_{2}\mu_{3}\alpha_{2}+\mu_{2}\mu_{3}\beta)a,
\end{align*}
Therefore,
\begin{align*}
[x,[y,z]]&=[\lambda_{1}a+\mu_{1}b,(\lambda_{2}\lambda_{3}\sigma+\mu_{2}\lambda_{3}\alpha_{1}+\lambda_{2}\mu_{3}\alpha_{2}+\mu_{2}\mu_{3}\beta)a]\\
&=\lambda_{1}(\lambda_{2}\lambda_{3}\sigma+\mu_{2}\lambda_{3}\alpha_{1}+\lambda_{2}\mu_{3}\alpha_{2}+\mu_{2}\mu_{3}\beta)[a,a]\\
&+\mu_{1}(\lambda_{2}\lambda_{3}\sigma+\mu_{2}\lambda_{3}\alpha_{1}+\lambda_{2}\mu_{3}\alpha_{2}+\mu_{2}\mu_{3}\beta)[b,a]\\
&=\lambda_{1}(\lambda_{2}\lambda_{3}\sigma+\mu_{2}\lambda_{3}\alpha_{1}+\lambda_{2}\mu_{3}\alpha_{2}+\mu_{2}\mu_{3}\beta)\sigma a\\
&+\mu_{1}(\lambda_{2}\lambda_{3}\sigma+\mu_{2}\lambda_{3}\alpha_{1}+\lambda_{2}\mu_{3}\alpha_{2}+\mu_{2}\mu_{3}\beta)\alpha_{1}a\\
&=(\mu_{1}\mu_{2}\lambda_{3}\alpha_{1}^{2}+\mu_{1}\lambda_{2}\mu_{3}\alpha_{2}\alpha_{1}+\mu_{1}\mu_{2}\mu_{3}\beta\alpha_{1})a,\\
[[x,y],z]&=[(\lambda_{1}\lambda_{2}\sigma+\mu_{1}\lambda_{2}\alpha_{1}+\lambda_{1}\mu_{2}\alpha_{2}+\mu_{1}\mu_{2}\beta)a,\lambda_{3}a+\mu_{3}b]\\
&=\lambda_{3}(\lambda_{1}\lambda_{2}\sigma+\mu_{1}\lambda_{2}\alpha_{1}+\lambda_{1}\mu_{2}\alpha_{2}+\mu_{1}\mu_{2}\beta)[a,a]\\
&+\mu_{3}(\lambda_{1}\lambda_{2}\sigma+\mu_{1}\lambda_{2}\alpha_{1}+\lambda_{1}\mu_{2}\alpha_{2}+\mu_{1}\mu_{2}\beta)[a,b]\\
&=\lambda_{3}(\lambda_{1}\lambda_{2}\sigma+\mu_{1}\lambda_{2}\alpha_{1}+\lambda_{1}\mu_{2}\alpha_{2}+\mu_{1}\mu_{2}\beta)\sigma a\\
&+\mu_{3}(\lambda_{1}\lambda_{2}\sigma+\mu_{1}\lambda_{2}\alpha_{1}+\lambda_{1}\mu_{2}\alpha_{2}+\mu_{1}\mu_{2}\beta)\alpha_{2}a\\
&=(\lambda_{3}\lambda_{1}\lambda_{2}\sigma^{2}+\lambda_{3}\mu_{1}\lambda_{2}\alpha_{1}\sigma+\lambda_{3}\lambda_{1}\mu_{2}\alpha_{2}\sigma+\lambda_{3}\mu_{1}\mu_{2}\beta\sigma\\
&+\mu_{3}\lambda_{1}\lambda_{2}\sigma\alpha_{2}+\mu_{3}\mu_{1}\lambda_{2}\alpha_{1}\alpha_{2}+\mu_{3}\alpha_{1}\mu_{2}\alpha_{2}^{2}+\mu_{3}\mu_{1}\mu_{2}\beta\alpha_{2})a\\
&=(\mu_{3}\mu_{1}\lambda_{2}\alpha_{1}\alpha_{2}+\mu_{3}\lambda_{1}\mu_{2}\alpha_{2}^{2})a,\\
[y,[x,z]]&=[\lambda_{2}a+\mu_{2}b,(\lambda_{1}\lambda_{3}\sigma+\mu_{1}\lambda_{3}\alpha_{1}+\lambda_{1}\mu_{3}\alpha_{2}+\mu_{1}\mu_{3}\beta)a]\\
&=\lambda_{2}(\lambda_{1}\lambda_{3}\sigma+\mu_{1}\lambda_{3}\alpha_{1}+\lambda_{1}\mu_{3}\alpha_{2}+\mu_{1}\mu_{3}\beta)[a,a]\\
&+\mu_{2}(\lambda_{1}\lambda_{3}\sigma+\mu_{1}\lambda_{3}\alpha_{1}+\lambda_{1}\mu_{3}\alpha_{2}+\mu_{1}\mu_{3}\beta)[b,a]\\
&=\lambda_{2}(\lambda_{1}\lambda_{3}\sigma+\mu_{1}\lambda_{3}\alpha_{1}+\lambda_{1}\mu_{3}\alpha_{2}+\mu_{1}\mu_{3}\beta)\sigma a\\
&+\mu_{2}(\lambda_{1}\lambda_{3}\sigma+\mu_{1}\lambda_{3}\alpha_{1}+\lambda_{1}\mu_{3}\alpha_{2}+\mu_{1}\mu_{3}\beta)\alpha_{1}a\\
&=(\lambda_{2}\lambda_{1}\lambda_{3}\sigma^{2}+\lambda_{2}\mu_{1}\lambda_{3}\alpha_{1}\sigma+\lambda_{2}\lambda_{1}\mu_{3}\alpha_{2}\sigma+\lambda_{2}\mu_{1}\mu_{3}\beta\sigma\\
&+\mu_{2}\lambda_{1}\lambda_{3}\sigma\alpha_{1}+\mu_{2}\mu_{1}\lambda_{3}\alpha_{1}^{2}+\mu_{2}\lambda_{1}\mu_{3}\alpha_{2}\alpha_{1}+\mu_{2}\mu_{1}\mu_{3}\beta\alpha_{1})a\\
&=(\mu_{2}\mu_{1}\lambda_{3}\alpha_{1}^{2}+\mu_{2}\lambda_{1}\mu_{3}\alpha_{2}\alpha_{1}+\mu_{2}\mu_{1}\mu_{3}\beta\alpha_{1})a.
\end{align*}

We have
\begin{gather*}
[x,[y,z]]-[[x,y],z]-[y,[x,z]]\\
=(\mu_{1}\mu_{2}\lambda_{3}\alpha_{1}^{2}+\mu_{1}\lambda_{2}\mu_{3}\alpha_{2}\alpha_{1}+\mu_{1}\mu_{2}\mu_{3}\beta\alpha_{1})a\\
-(\mu_{3}\mu_{1}\lambda_{2}\alpha_{1}\alpha_{2}+\mu_{3}\lambda_{1}\mu_{2}\alpha_{2}^{2})a\\
-(\mu_{2}\mu_{1}\lambda_{3}\alpha_{1}^{2}+\mu_{2}\lambda_{1}\mu_{3}\alpha_{2}\alpha_{1}+\mu_{2}\mu_{1}\mu_{3}\beta\alpha_{1})a\\
=-(\mu_{3}\lambda_{1}\mu_{2}\alpha_{2}^{2}+\mu_{2}\lambda_{1}\mu_{3}\alpha_{2}\alpha_{1})a\\
=-\mu_{3}\lambda_{1}\mu_{2}(\alpha_{2}^{2}+\alpha_{2}\alpha_{1})a=0.
\end{gather*}
Thus we obtain
$$[x,[y,z]]=[[x,y],z]+[y,[x,z]],$$
which shows that $L$ is a Leibniz ring.
By the above, we have
$$[x,y]=(\lambda_{1}\lambda_{2}\sigma+\mu_{1}\lambda_{2}\alpha_{1}+\lambda_{1}\mu_{2}\alpha_{2}+\mu_{1}\mu_{2}\beta)a.$$
Hence
$$[L,L]=\langle\sigma a\rangle+\langle\alpha_{1}a\rangle+\langle\alpha_{2}a\rangle+\langle\beta a\rangle.$$
An element $x\in\zeta^{\mathrm{left}}(L)$ if and only if $[x,a]=[x,b]=0$. We have
\begin{align*}
[x,a]&=(\lambda_{1}\sigma+\mu_{1}\alpha_{1})a,\\
[x,b]&=(\lambda_{1}\alpha_{2}+\mu_{1}\beta)a.
\end{align*}
It follows that
$$\zeta^{\mathrm{left}}(L)=\langle\sigma a\rangle+\langle\alpha_{1}a\rangle+\langle\alpha_{2}a\rangle+\langle\beta a\rangle=[L,L].$$

Similarly, $y\in\zeta^{\mathrm{right}}(L)$ if and only if $[a,y]=[b,y]=0$. We have
\begin{align*}
[a,y]&=(\lambda_{2}\sigma+\mu_{2}\alpha_{2})a,\\
[y,b]&=(\lambda_{2}\alpha_{1}+\mu_{2}\beta)a,
\end{align*}
and again conclude that
$$\zeta^{\mathrm{right}}(L)=\langle\sigma a\rangle+\langle\alpha_{1}a\rangle+\langle\alpha_{2}a\rangle+\langle\beta a\rangle=[L,L]$$
Thus,
$$\zeta^{\mathrm{left}}(L)=\zeta^{\mathrm{right}}(L)=\zeta(L)=[L,L].$$
Furthermore,
$$[x,x]=(\lambda_{1}^{2}\sigma+\mu_{1}\lambda_{1}\alpha_{1}+\lambda_{1}\mu_{1}\alpha_{2}+\mu_{1}^{2}\beta)a,$$
and hence
$$K=\langle\sigma a\rangle+\langle(\alpha_{1}+\alpha_{2})a\rangle+\langle\beta a\rangle.$$
Therefore, we obtain the following type of Leibniz ring:
\begin{gather*}
L_{9}=\langle a_{1}\rangle\oplus\langle a_{2}\rangle,\\
\mbox{the additive order of}\ a_{1}\ \mbox{is finite},\ |a_{1}|=k,\\
\mbox{the additive order of}\ a_{2}\ \mbox{is infinite},\\
[a_{1},a_{1}]=\sigma a_{1},\ [a_{1},a_{2}]=\alpha_{2}a,\ [a_{2},a_{1}]=\alpha_{1}a_{1},\ [a_{2},a_{2}]=\beta a_{1},\\
\mbox{where}\ \sigma,\alpha_{1},\alpha_{2},\beta\ \mbox{are integers, satisfying the following conditions:}\\
\sigma^{2}\equiv0(\mathrm{mod}\ k),\ \alpha_{2}\sigma\equiv0(\mathrm{mod}\ k),\alpha_{2}\beta\equiv0(\mathrm{mod}\ k),\\ \beta\sigma\equiv0(\mathrm{mod}\ k),\alpha_{2}^{2}+\alpha_{2}\alpha_{1}\equiv0(\mathrm{mod}\ k),\ \alpha_{1}\sigma\equiv0(\mathrm{mod}\ k).
\end{gather*}
Here
\begin{gather*}
\zeta(L_{9})=\zeta^{\mathrm{left}}(L_{9})=\zeta^{\mathrm{right}}(L_{9})=[L_{9},L_{9}]\\
=\langle\sigma a\rangle+\langle\alpha_{1}a\rangle+\langle\alpha_{2}a\rangle+\langle\beta a\rangle,\\
Leib(L_{9})=\langle\sigma a\rangle+\langle(\alpha_{1}+\alpha_{2})a\rangle+\langle\beta a\rangle.
\end{gather*}


\begin{thebibliography}{OO}
\bibitem{AOR2020} Ayupov, Sh., Omirov, B., Rakhimov, I.: Leibniz Algebras: Structure and Classification. CRC Press, Taylor \& Francis Group (2020). https://doi.org/10.1201/9780429344336
\bibitem{BD2013} Barnes, D.: Schunck classes of soluble Leibniz algebras. Comm. Algebra \textbf{41}(11), 4046-4065 (2013). https://doi.org/10.1080/00927872. 2012.700978
\bibitem{BA1965} Blokh, A.: On a generalization of the concept of Lie algebra. Dokl. Akad. Nauk SSSR \textbf{165}(3), 471-473 (1965).
\bibitem{BA1967} Blokh, A.: Cartan-Eilenberg homology theory for a generalized class of Lie algebras. Dokl. Akad. Nauk SSSR \textbf{175}(8), 824-826 (1967).
\bibitem{BA1971} Blokh, A.: A certain generalization of the concept of Lie algebra. Algebra and number theory. Moskov. Gos. Ped. Inst. Uchen. Zap. \textbf{375}, 9-20 (1971).
\bibitem{CPSeY2019} Chupordia, V.A., Pypka, A.A., Semko, N.N., Yashchuk, V.S.: Leibniz algebras: a brief review of current results. Carpathian Math. Publ. \textbf{11}(2), 250-257 (2019). https://doi.org/10.15330/cmp.11.2.250-257
\bibitem{KiKuPSu2017} Kirichenko, V.V., Kurdachenko, L.A., Pypka, A.A., Subbotin, I.Ya.: Some aspects of Leibniz algebra theory. Algebra Discrete Math. \textbf{24}(1), 1-33 (2017).
\bibitem{KOP2016} Kurdachenko, L.A., Otal, J., Pypka, A.A.: Relationships between the factors of the canonical central series of Leibniz algebras. Eur. J. Math. \textbf{2}(2), 565-577 (2016). https://doi.org/10.1007/s40879-016-0093-5
\bibitem{KPS2018} Kurdachenko, L.A., Pypka, A.A., Subbotin, I.Ya.: On new analogs of some classical group theoretical results in Lie rings. In ``Infinite Group Theory: From the Past to the Future'', Chapter 11, 197-213 (2018). https://doi.org/10.1142/9789813204058\_0011
\bibitem{KPS2021} Kurdachenko, L.A., Pypka, O.O., Subbotin, I.Ya.: On Leibniz algebras whose subalgebras are either ideals or self-idealizing subalgebras. Ukrainian Math. J. \textbf{73}(6), 944-962 (2021). https://doi.org/10.1007/ s11253-021-01969-0
\bibitem{KPS2022} Kurdachenko, L.A., Pypka, O.O., Subbotin, I.Ya.: On the structure of low-dimensional Leibniz algebras: some revision. Algebra Discrete Math. \textbf{34}(1), 68-104 (2022). https://doi.org/10.12958/adm2036
\bibitem{KPS2023} Kurdachenko, L.A., Pypka, A.A., Subbotin, I.Ya.: On the automorphism groups of some Leibniz algebras. Int. J. Group Theory. \textbf{12}(1), 1-20 (2023). https://doi.org/10.22108/IJGT.2021.130057.1735
\bibitem{KPS2024} Kurdachenko, L.A., Pypka, O.O., Subbotin, I.Ya.: General Theory of Leibniz Algebras. Switzerland: Springer (2024). https://doi.org/10. 1007/978-3-031-58148-9
\bibitem{KPV2023} Kurdachenko, L.A., Pypka, O.O., Velychko, T.V.: On the automorphism groups for some Leibniz algebras of low dimensions. Ukrainian Math. J. \textbf{74}(10), 1526-1546 (2023). https://doi.org/10.1007/s11253-023-02153-2
\bibitem{KSS2017} Kurdachenko, L.A., Semko, N.N., Subbotin, I.Ya.: The Leibniz algebras whose subalgebras are ideals. Open Math. \textbf{15}(1), 92-100 (2017). https://doi.org/10.1515/math-2017-0010
\bibitem{KSS2018A} Kurdachenko, L.A., Semko, N.N., Subbotin, I.Ya.: On the anticommutativity in Leibniz algebras. Algebra Discrete Math. \textbf{26}(1), 97--109 (2018).
\bibitem{KSeSu2020} Kurdachenko, L.A., Semko, N.N., Subbotin, I.Ya.: Applying group theory philosophy to Leibniz algebras: some new developments. Adv. Group Theory Appl. \textbf{9}, 71-121 (2020). https://doi.org/10.32037/agta-2020-004
\bibitem{KSeSu2023} Kurdachenko, L.A., Semko, M.M., Subbotin, I.Ya.: On the algebra of derivations of some low-dimensional Leibniz algebras. Algebra Discrete Math. \textbf{36}(1), 43-60 (2023). https://doi.org/10.12958/adm2161
\bibitem{KSeSu2024} Kurdachenko, L.A., Semko, M.M., Subbotin, I.Ya.: On the structure of algebras of derivations of some non-nilpotent Leibniz algebras. Algebra Discrete Math. \textbf{37}(2), 244-261 (2024). https://doi.org/10.12958/ adm2227
\bibitem{KSeSu2024A} Kurdachenko, L.A., Semko, M.M., Subbotin, I.Ya.: On the algebra of derivations of some Leibniz algebras. Algebra Discrete Math. \textbf{38}(1), 63-86 (2024). https://doi.org/10.12958/adm2316
\bibitem{KSeY2021} Kurdachenko, L.A., Semko, M.M., Yashchuk, V.S.: On the structure of the algebra of derivations of cyclic Leibniz algebras. Algebra Discrete Math. \textbf{32}(2), 241-252 (2021). https://doi.org/10.12958/adm1898
\bibitem{KSeYa2023} Kurdachenko, L.A., Semko, M.M., Yashchuk, V.S.: On the algebra of derivations of some nilpotent Leibniz algebras. Res. Math. \textbf{31}(1), 62-71 (2023). https://doi.org/10.15421/242306
\bibitem{KSeYa2024} Kurdachenko, L.A., Semko, M.M., Yashchuk, V.S.: On the structure of the algebra of derivations of some low-dimensional Leibniz algebras. Ukrainian Math. J. \textbf{76}(5), 728-742 (2024). https://doi.org/10.3842/ umzh.v76i5.7573
\bibitem{KSuSe2018} Kurdachenko, L.A., Subbotin, I.Ya., Semko, N.N.: From groups to Leibniz algebras: Common approaches, parallel results. Adv. Group Theory Appl. \textbf{5}, 1-31 (2018). https://doi.org/10.4399/97888255161421
\bibitem{KSuY2024} Kurdachenko, L.A., Subbotin, I.Ya., Yashchuk, V.S.: On the endomorphisms and derivations of some Leibniz algebras. J. Algebra Appl. \textbf{23}(1), 2450002 (2024). https://doi.org/10.1142/S0219498824500026
\bibitem{LJ1992} Loday, J.-L.: Cyclic homology. Grundlehren der Mathematischen Wissenschaften. \textbf{301}. Springer Verlag (1992). https://doi.org/10.1007/ 978-3-662-11389-9
\bibitem{LJ1993} Loday, J.-L.: Une version non commutative des alg\`{e}bres de Lie; les alg\`{e}bras de Leibniz. Enseign. Math. \textbf{39}, 269-293 (1993).
\bibitem{P2022} Pypka, O.O.: Some relationships between the generalized central series of Leibniz algebras. Ukrainian Math. J. \textbf{73}(12), 1958-1966 (2022). https://doi.org/10.1007/s11253-022-02040-2
\bibitem{SSY2022} Semko, M.M., Skaskiv, L.V., Yarovaya, O.A.: On the derivations of cyclic nilpotent Leibniz algebras. Carpathian Math. Publ. \textbf{14}(2), 345-353 (2022). https://doi.org/10.15330/cmp.14.2.345-353.
\end{thebibliography}
\end{document}